\documentclass{amsart}

\DeclareMathSymbol{\twoheadrightarrow} {\mathrel}{AMSa}{"10}

\def\RR{{\mathfrak R}}
\def\F{{\mathbf F}}
\def\L{{\mathbf L}}

\def\SS{{\mathbf S}}
\def\A{{\mathbf A}}
\def\AGL{\mathrm{AGL}}
\def\UU{\mathrm{U}}
\def\Sn{{\mathbf S}_n}

\def\Gal{\mathrm{Gal}}
\def\Perm{\mathrm{Perm}}

\def\Pr{\mathrm{Pr}}

\def\End{\mathrm{End}}
\def\Aut{\mathrm{Aut}}

\def\Mat{\mathrm{Mat}}

\def\I{\mathrm{Id}}

\def\fchar{\mathrm{char}}

\def\GL{\mathrm{GL}}
\def\PGL{\mathrm{PGL}}
\def\PSL{\mathrm{PSL}}

\def\SU{\mathrm{SU}}

\def\M{\mathrm{M}}

\def\dim{\mathrm{dim}}

\def\P{{\mathbf P}}
\def\Sz{\mathrm{Sz}}

\newtheorem{thm}{Theorem}[section]
\newtheorem{lem}[thm]{Lemma}
\newtheorem{cor}[thm]{Corollary}

\theoremstyle{definition}
\newtheorem{defn}[thm]{Definition}
\newtheorem{ex}[thm]{Example}
\newtheorem{exs}[thm]{Examples}

\newtheorem{rem}[thm]{Remark}
\newtheorem{rems}[thm]{Remarks}
\title[Very simple representations]
{Very simple representations: variations on a theme of Clifford}
\author[Yuri G. Zarhin]{Yuri G. Zarhin}
\address{Department of Mathematics, Pennsylvania State University,
University Park, PA 16802, USA} \email{zarhin\char`\@math.psu.edu}
\thanks{Partially supported by the NSF}

\begin{document}

            \begin{abstract}
  We discuss a certain class of absolutely irreducible group
  representations that  behave nicely under the restriction
  to normal subgroups and  subalgebras. Interrelations with doubly
  transitive permutation groups and hyperelliptic jacobians are
  discussed.
 \end{abstract}

            \maketitle

\section{Introduction}

We start this paper with the following natural definition.

\begin{defn}
Let $V$ be a vector space over a field $k$, let $G$ be a group and
$\rho: G \to \Aut_{k}(V)$ a linear representation of $G$ in $V$.
Suppose $R\subset \End_{k}(V)$ is an $k$-subalgebra containing the
identity operator $\I$. We say that $R$ is {\sl normal} if
$$\rho(s) R \rho(s)^{-1} \subset R \quad \forall s \in G.$$
\end{defn}

\begin{exs}
Clearly, $\End_{k}(V)$ is normal. The algebra $k\cdot \I$ of
scalars is also normal.
 If $H$ is a {\sl normal} subgroup of $G$ then the image of the group
algebra $k[H]$ in $\End_{k}(V)$ is normal.
\end{exs}

The following assertion is a straightforward generalization of
well-known Clifford's theorem \cite{Clifford}; ~\cite[\S 49]{CR};
~\cite[\S 8.1]{Serre}; ~\cite[Ch. 6]{Isaacs}.

\begin{lem}[Lemma 7.4 of \cite{ZarhinTexel}]
 \label{induce} Let $G$ be a group, $k$ a field, $V$ a non-zero $k$-vector
space of finite dimension $n$ and
$$\rho: G \to \Aut_k(V)$$
an irreducible representation.
  Let $R \subset \End_{\F}(V)$ be a normal subalgebra.

Then:
\begin{enumerate}
\item[(i)]
The faithful $R$-module $V$ is semisimple.
\item[(ii)]
Either the $R$-module $V$ is isotypic or there exists a subgroup
$G' \subset G$ of index $r$ dividing $n$ and a $G'$-module $V'$ of
finite $k$-dimension $n/r$ such that
    $r>1$ and
the $G$-module $V$ is induced from $V'$.
\end{enumerate}
\end{lem}

The following notion was introduced by the author in
\cite{ZarhinTexel} (see also \cite{ZarhinMRL2}); it  proved to be
useful for the construction of abelian varieties with small
endomorphism rings \cite{ZarhinMMJ,ZarhinCrelle}.

\begin{defn}
Let $V$ be a non-zero finite-dimensional vector space over a field
$k$, let $G$ be a group and $\rho: G \to \Aut_{k}(V)$ a linear
representation of $G$ in $V$. We say that the $G$-module $V$ is
{\sl very simple} if it enjoys the following property:

If $R \subset \End_{k}(V)$ is a normal subalgebra
 then either $R=k\cdot \I$ or $R=\End_{k}(V)$.

 Here is (obviously) an equivalent definition: if $R \subset \End_{k}(V)$ is a normal subalgebra
  then either $\dim_k(R)=1$ or $\dim_k(R)=(\dim_k(V))^2$.
\end{defn}

In this paper we prove that very simple representations over an
algebraically closed field are exactly those absolutely
irreducible representations that are {\sl not} induced from a
representation of a proper subgroup and do {\sl not split}
non-trivially into a tensor product of projective representations.
This assertion remains valid  for representations of {\sl perfect}
groups over {\sl finite} fields. We also give a certain criterion
that works for any ground field with trivial Brauer group.

The paper is organized as follows. In \S \ref{basic} we list basic
properties of very simple representations. In \S \ref{counter} we
discuss certain natural constructions of representations that are
{\sl not} very simple. The Section \ref{main} contains the
statement of main results about very simple representations and
their proof.
%%%%%%%%%%%%%%%%%%%%%%%%%%%%%%%%%%%%%%%%%%%%%%%%%%%%%%%%%%%%%%%%%%%
In the Section \ref{permute} we discuss interrelations between
very simple representations and doubly transitive permutation
groups. The last section contains applications to hyperelliptic
jacobians.
%%%%%%%%%%%%%%%%%%%%%%%%%%%%%%%%%%%%%%%%%%%%%%%%%%%%%%%%%%%%%%%%%%%%%

{\bf Acknowledgements}. This paper germinated during author's
short stay in Manchester in August of 2002. The author would like
to thank UMIST  Department of Mathematics for its hospitality. My
special thanks go to Professor A. V. Borovik for helpful
encouraging discussions.

\section{Very simple representations}
\label{basic}

\begin{rems}
\label{image}
\begin{enumerate}
\item[(i)]
Clearly, if $\dim_k(V)=1$ then the $G$-module $V$ is very simple.
In other words, every one-dimensional representation is very
simple.
\item[(ii)]
Clearly, the $G$-module $V$ is very simple if and only if the
corresponding $\rho(G)$-module $V$ is very simple.
\item[(iii)]
Clearly, if $V$ is very simple then the corresponding algebra
homomorphism
    $$k[G] \to \End_{k}(V)$$
is surjective. Here $k[G]$ stands for the group algebra of $G$. In
particular, a very simple module is absolutely simple.
\item[(iv)]
If $G'$ is a subgroup of $G$ and the $G'$-module $V$ is very
simple then the $G$-module $V$ is also very simple.
\item[(v)] Suppose $W$ is a one-dimensional $k$-vector space and
$$\kappa: G \to k^*=\Aut_k(W)$$
is a one-dimensional representation of $G$. Then the $G$-module
$V\otimes_k W$ is very simple if and only if the $G$-module $V$ is
very simple. Indeed, there are the canonical $k$-algebra
isomorphisms
$$\End_k(V)=\End_k(V)\otimes_k k=\End_k(V)\otimes_k \End_k(W) \cong \End_k(V\otimes_k W),$$
which  are isomorphisms of the corresponding $G$-modules.

\item[(vi)]
Let $G'$ be a normal subgroup of $G$. If $V$ is a very simple
$G$-module then either $\rho(G') \subset \Aut_{k}(V)$ consists of
scalars (i.e., lies in $k\cdot\I$) or the $G'$-module $V$ is
absolutely simple. Indeed, let $R'\subset \End_{k}(V)$ be the
image of the natural homomorphism $k[G'] \to \End_{k}(V)$.
Clearly, $R'$ is normal. Hence either $R'$ consists of scalars and
therefore $\rho(G')\subset R'$ consists of scalars or $R'=
\End_{k}(V)$ and therefore the $G'$-module $V$ is absolutely
simple.

As an immediate corollary, we get the following assertion: if
$\dim_{k}(V)>1$, the $G$-module $V$ is {\sl faithful} very simple
and $G'$ is a {\sl non}-central normal subgroup of $G$ then the
$G'$-module $V$ is absolutely simple; in particular, $G'$ is {\sl
non}-abelian. In addition, if $k=\F_2$ then the only abelian
normal subgroup of $G$ is the trivial (one-element) subgroup.
\item[(vii)]
Suppose $F$ is a discrete valuation field with  valuation ring
$O_F$, maximal ideal $m_F$ and residue field $k=O_F/m_F$. Suppose
$V_F$ is a finite-dimensional $F$-vector space and
$$\rho_F: G \to \Aut_F(V_F)$$
is a $F$-linear representation of $G$. Suppose $T$ is a $G$-stable
$O_F$-lattice in $V_F$ and the corresponding $k[G]$-module $T/m_F
T$ is isomorphic to $V$. Assume that the $G$-module $V$ is very
simple. Then the $G$-module $V_F$ is also very simple. In other
words, a lifting of a very simple module is also very simple. (See
\cite{ZarhinMMJ}, Remark 5.2(v).)
\end{enumerate}
\end{rems}

\begin{ex}
\label{dim22}
 Suppose $k=\F_2, \dim_k(V)=2, G=\Aut_k(V)=\GL_2(\F_2)$. Then the faithful absolutely simple $G$-module $V$ is not
 very simple. Indeed, $G$ is isomorphic to the symmetric group
 $\SS_3$ and therefore contains a non-central abelian normal subgroup
 isomorphic to the alternating group $\A_3$. By Remark
 \ref{image}(vi), the $G$-module $V$ is not very simple.

 Applying Remarks \ref{image}(ii) and \ref{image}(iv), we conclude
 that all two-dimensional representations over $\F_2$ of any group
 are not very simple.
\end{ex}

\begin{ex}
\label{small}
 Suppose $V$ is a finite-dimensional vector space
over a {\sl finite} field $k$ of characteristic $\ell$ and $G$ is
a {\sl perfect} subgroup of $\Aut(V)$ enjoying the following
properties:
\begin{enumerate}
\item[(i)]
If $Z$ is the center of $G$  then the quotient $\Gamma:=G/Z$ is a
simple non-abelian group.
\item[(ii)]
Every nontrivial projective representation of $\Gamma$ in
characteristic $\ell$ has dimension $\ge \dim_k(V)$.
\end{enumerate}
Then  the $G$-module $V$ is very simple. See Cor. 5.4 in
\cite{ZarhinMMJ}.
\end{ex}

\section{Counterexamples}
\label{counter}
 Throughout this section, $k$ is a field, $V$ a non-zero
 finite-dimensional
$k$-vector space and $\rho: G \to \Aut_k(V)$ is a linear
representation of $G$ in $V$.

\begin{ex}
\begin{enumerate}
\item[(i)]
Assume that there exist $k[G]$-modules $V_1$ and $V_2$ such that
$\dim_k(V_1)>1, \dim_k(V_2)>1$ and the $G$-module $V$ is
isomorphic to $V_1 \otimes_{k} V_2$. Then $V$ is {\sl not} very
simple. Indeed, the subalgebra
$$R=\End_{k}(V_1)\otimes  \I_{V_2}\subset
\End_{k}(V_1)\otimes_{k}\End_{k}(V_2)= \End_{k}(V_1 \otimes_{k}
V_2)=\End_{k}(V)$$ is normal
 but coincides neither
with $k\cdot \I$ nor with $\End_{k}(V)$. (Here $\I_{V_2}$ stands
for the identity operator in $V_2$.) Clearly, the centralizer of
$R$ in $\End_{k}(V)$ coincides with $\I_{V_1}\otimes\End_{k}(V_2)$
and is also normal. (Here $\I_{V_1}$ stands for the identity
operator in $V_1$.)

\item[(ii)]
Let $X \twoheadrightarrow G$ be a {\sl surjective} group
homomorphism. Assume that there exist $k[X]$-modules $V_1$ and
$V_2$ such that $\dim_k(V_1)>1, \dim_k(V_2)>1$ and $V$, viewed as
$X$-module, is isomorphic to $V_1 \otimes_{k} V_2$. Then the
$X$-module $V$ is {\sl not} very simple. Since $X$ and $G$ have
the same image in $\Aut_{k}(V)$, the $G$-module $V$ is also {\sl
not} very simple.
\end{enumerate}
\end{ex}

\begin{defn}
\label{split} Let $V$ be a vector space over a field $k$, let $G$
be a group and $\rho: G \to \Aut_{k}(V)$ a linear representation
of $G$ in $V$. Let
$$1 \to C \hookrightarrow X \stackrel{\pi}{\twoheadrightarrow} G
\to 1$$ be a {\sl central} extension of $G$, i.e., $C$ is a {\sl
central} subgroup of $G$ which coincides with the kernel of {\sl
surjective} homomorphism $\pi:X \to G$.
 Suppose that the representation
$$X \stackrel{\pi}{\twoheadrightarrow} G \stackrel{\rho}{\longrightarrow}  \Aut_k(V)$$ of $X$
is isomorphic to a tensor product $\rho_1\otimes \rho_2:X \to
\Aut_k(V_1\otimes_k V_2)$ of two $k$-linear representations
$$\rho_1:X \to \Aut_k(V_1), \quad \rho_2:X \to \Aut_k(V_2)$$
with $$\dim_k(V_1)>1, \quad \dim_k(V_2)>1.$$ Then we say that the
$G$-module $V$ {\sl splits} and call the triple
$(X\stackrel{\pi}{\twoheadrightarrow} G; \rho_1, \rho_2)$ a {\sl
splitting} of the $G$-module $V$.

We say that $V$ splits {\sl projectively} if both images
$\rho_1(C)\subset\Aut_k(V_1)$ and $\rho_2(C)\subset\Aut_k(V_2)$
consist of scalars. In other words, one may view both $\rho_1$ and
$\rho_2$ as {\sl projective} representations of $G$. In this case
we call  $(X\stackrel{\pi}{\twoheadrightarrow} G; \rho_1, \rho_2)$
a {\sl projective splitting} of the $G$-module $V$.

We call a splitting $(X\stackrel{\pi}{\twoheadrightarrow} G;
\rho_1, \rho_2)$ {\sl absolutely simple} if both $\rho_1$ and
$\rho_2$ are absolutely irreducible representations of $X$.
\end{defn}

\begin{rems}
\label{proj} We keep the notations of Definition \ref{split}.
\begin{enumerate}
\item[(i)]
Clearly,
$$\End_X(V_1)\otimes_k\End_X(V_2)\subset \End_X(V_1\otimes_k
V_2)=\End_X(V)=\End_G(V).$$ This implies that if $\End_G(V)=k$
then $\End_X(V_1)=k$ and $\End_X(V_2)=k$.
\item[(ii)]
Suppose $W$ is a {\sl proper} $X$-invariant subspace in $V_1$
(resp. in $V_2$). Then $W\otimes_k V_2$ (resp. $V_1\otimes_k W$)
is a {\sl proper} $X$-invariant subspace in $V_1\otimes_k V_2=V$
and therefore the corresponding $X$-module $V$ is {\sl not}
simple. This implies that the $G$-module $V$ is also {\sl not}
simple.

\item[(iii)]
 Suppose that the $G$-module $V$ is  absolutely simple
and splits. It follows from (i) and (ii) that both $\rho_1$ and
$\rho_2$ are also absolutely simple. In other words, every
splitting of an absolutely simple module is absolutely simple.

Now the centrality of $C$ combined with the absolute
irreducibility of $\rho_1$ and $\rho_2$ implies, thanks to Schur's
Lemma, that both images $\rho_1(C)\subset\Aut_k(V_1)$ and
$\rho_2(C)\subset\Aut_k(V_2)$ consist of scalars. In other words,
every splitting of an absolutely simple $G$-module is projective.

\item[(iv)]
Suppose $G$ is a finite {\sl perfect} group. We write
$\gamma:\tilde{G}\twoheadrightarrow G$ for the universal central
extension of $G$~\cite[Ch. 2, \S 9]{Suzuki} also known as the {\sl
representation group} or the {\sl primitive central extension} of
$G$. It is known~\cite[Ch. 2, Th. 9.18]{Suzuki} that  $\tilde{G}$
is also a finite perfect group and for each central extension
$X\stackrel{\pi}{\twoheadrightarrow} G$  there exists a surjective
homomorphism $\phi:\tilde{G}\twoheadrightarrow [X,X]$ to the
derived subgroup $[X,X]$ of $X$ such that the composition
$$\tilde{G}\stackrel{\phi}{\twoheadrightarrow} [X,X]\subset
X \stackrel{\pi}{\twoheadrightarrow} G$$ coincides with
$\gamma:\tilde{G}\twoheadrightarrow G$. This implies that while
checking whether the $G$-module $V$ admits a projective absolutely
simple splitting,  one may always restrict oneself to the case of
$X=\tilde{G}$ and $\pi=\gamma$, and deal exclusively with
absolutely irreducible linear representations of $\tilde{G}$ over
$k$.
\item[(v)]
We refer the reader to \cite{Rachel} for a study of projective
representations of arbitrary finite groups over not necessarily
algebraically closed fields.
\end{enumerate}
\end{rems}

\begin{ex}
\label{ind}
 Assume that there exists a subgroup $G'$ in $G$ of
finite index $m>1$ and a $G'$-module $W$ such that the
$k[G]$-module $V$ is {\sl induced} from the $k[G']$-module $W$.
Then $V$ is {\sl not} very simple. Indeed, one may view $W$ as a
$G'$-submodule of $V$ such that $V$ coincides with the direct sum
$\oplus_{\sigma \in G/G'}\sigma W$ and $G$ permutes all $\sigma
W$'s. We write $\Pr_{\sigma}:V\twoheadrightarrow \sigma W\subset
V$ for the corresponding projection maps. Then
 $R=\oplus_{\sigma\in G/G'}k\cdot\Pr_{\sigma}$ is the algebra of
all operators sending each $\sigma W$ into itself and acting on
each $\sigma W$ as scalars.  Clearly, $R$ is normal but coincides
neither with $k\cdot \I$ nor with $\End_{k}(V)$.

Notice that if the $G'$-module $V'$ is {\sl trivial} (i.e.,
$s(w)=w$ for all $s\in G', w\in W$)  then the $G$-module $V$ is
{\sl not} simple. Indeed, for any non-zero $w\in W$ the vector
$$v=\sum_{\sigma \in G/G'}\sigma(w) \in \oplus_{\sigma \in
G/G'}\sigma W=V$$ is a non-zero $G$-invariant element of $V$.
Since $\dim_k(V)\ge m>1$, the $G$-module $V$ is not simple.

Clearly, if $G'=\{1\}$ is the trivial subgroup of $G$ then every
$G'$ -module is trivial. This implies that if the conditions of
Lemma \ref{induce} hold true then (in the notations of
\ref{induce})
 either the $R$-module $V$ is  isotypic or the $G$-module
$V$ is induced from a representation of a {\sl proper} subgroup of
finite index.
\end{ex}

\begin{rem}
\label{semi}
 Suppose $k'/k$ is a finite algebraic extension of
fields. We write $\Aut(k'/k)$ for the group of $k$-linear
automorphisms of the field $k'$. It is well-known that
$\Aut(k'/k)$ is finite and its order divides the degree $[k':k]$;
the equality holds if and only if $k'/k$ is Galois.

Suppose there exists a homomorphism $\chi:G \to \Aut(k'/k)$
enjoying the following property:

There exists a structure of $k'$-vector space on $V$ such that
$$\rho(s)(av)=(\chi(s)(a))v \quad \forall s \in G, a \in k', v
\in V.$$ We claim that if the $G$-module $V$ is absolutely simple
then $k'/k$ is Galois and $\chi$ is surjective. Indeed, let us
consider the image $\chi(G)\subset \Aut(k'/k)$. We write $k_0$ for
the subfield of $\chi(G)$-invariants in $k'$. We have
$$k\subset k_0 \subset k';$$
the degree $[k':k_0]$ coincides with the order of $\chi(G)$ and
therefore divides the order of $\Aut(k'/k)$. Clearly, $G$ acts on
$V$ by $k_0$-linear automorphisms, i.e., $k_0$ commutes with the
action of $G$ on $V$. Now the absolute irreducibility of $V$
implies that $k_0=k$. This implies that $[k_0:k]=1$. Since
$[k':k]=[k':k_0][k_0:k]$, we conclude that $[k':k]=[k':k_0]$
coincides with the order of $\chi(G)$ and therefore divides the
order of $\Aut(k'/k)$. This implies that $[k':k]$ coincides with
the order of $\Aut(k'/k)$ and therefore $k'/k$ is Galois. Since
$[k':k]$ coincides with the order of $\chi(G)$, the order of the
group $\Aut(k'/k)$ must coincide with the order of its subgroup
$\chi(G)$ and therefore $\chi(G)= \Aut(k'/k)$, i.e., $\chi$ is
surjective.
\end{rem}

\begin{defn}
We say that the $G$-module $V$ admits a {\sl twisted
multiplication} if there exist a {\sl nontrivial} Galois extension
$k'$ of $k$ and a surjective homomorphism $\chi:G \to \Gal(k'/k)$
enjoying the following property:

There exists a structure of $k'$-vector space on $V$ such that
$$\rho(s)(av)=(\chi(s)(a))v \quad \forall s \in G, a \in k', v
\in V.$$ (Here $\Gal(k'/k)$ stands for the Galois group of
$k'/k$.) In other words, $G$ acts on $V$ by $k'$-semi-linear
automorphisms.
\end{defn}

\begin{ex}
Let us assume that $V$ admits a twisted multiplication. Then the
degree $[k':k]$ divides $\dim_k(V)$ and therefore $\dim_k(V)>1$.
Then the $G$-module $V$ is {\sl not} very simple. Indeed, $k'$ is
obviously normal but does coincide neither with $k\cdot \I$ nor
with $\End_k(V)$, since $k'\ne k$ and $\End_k(V)$ is
noncommutative.
\end{ex}

\begin{rem}
\label{nomult}
 Let us assume that  either $k$ is  algebraically
closed or $G$ is perfect and $k$ is finite. Then $V$ never admits
a twisted multiplication, because either every algebraic extension
$k'/k$ is trivial or $G$ is perfect and  every Galois group
$\Gal(k'/k)$ is abelian. In the latter case  every homomorphism
from perfect $G$ to abelian $\Gal(k'/k)$ must be trivial.
\end{rem}

\section{Main Theorem}
\label{main}

\begin{thm}
\label{central} Suppose the Brauer group of a field $k$ is trivial
(e.g., $k$ is either finite or algebraically closed). Suppose $V$
is a non-zero finite-dimensional $k$-vector space and
$$\rho:G \to \Aut_k(V)$$
is a linear representation of a group $G$ over $k$. Then the
$G$-module $V$ is  very simple if and only if all the following
conditions hold:
\begin{enumerate}
\item[(i)]
The $G$-module  $V$ is  absolutely simple;
\item[(ii)]
The $G$-module $V$ does not admit a projective absolutely simple
splitting;
\item[(iii)]
The $G$-module $V$ is not induced from a representation of a
proper subgroup  of finite index in $G$;
\item[(iv)]
The $G$-module $V$ does not admit a twisted multiplication.
\end{enumerate}
\end{thm}

\begin{cor}
\label{centralF} Let us assume that  either $k$ is  algebraically
closed or $G$ is perfect and $k$ is finite. Suppose $V$ is a
non-zero finite-dimensional $k$-vector space and
$$\rho:G \to \Aut_k(V)$$
is a linear representation of a group $G$ over $k$. Then the
$G$-module $V$ is  very simple if and only if all the following
conditions hold:
\begin{enumerate}
\item[(i)]
The $G$-module  $V$ is  absolutely simple;
\item[(ii)]
The $G$-module $V$ does not admit a projective absolutely simple
splitting;
\item[(iii)]
The $G$-module $V$ is not induced from a representation of a
proper subgroup of finite index in $G$.
\end{enumerate}
\end{cor}

\begin{proof}[Proof of Corollary \ref{centralF}]
Indeed, the Brauer group of $k$ is trivial. Now the proof follows
readily from Theorem \ref{central} combined with Remark
\ref{nomult}.

\end{proof}

Taking into account that every projective representation over
$\F_2$ is, in fact, linear, we obtain the following assertion.

\begin{cor}
\label{centralF2}  Suppose $V$ is a non-zero finite-dimensional
vector space over $\F_2$ and
$$\rho:G \to \Aut_{\F_2}(V)$$
is a linear representation of a group $G$ over $\F_2$. Then the
$G$-module $V$ is  very simple if and only if all the following
conditions hold:
\begin{enumerate}
\item[(i)]
The $G$-module  $V$ is  absolutely simple;
\item[(ii)]
The $G$-module $V$ does not split into a tensor product  $$V \cong
V_1\otimes_{\F_2}V$$ of two absolutely  simple $G$-modules $V_1$
and $V_2$ with $$\dim_{\F_2}(V_1)>1, \quad \dim_{\F_2}(V_2)>1;$$
\item[(iii)]
The $G$-module $V$ is not induced from a representation of a
proper subgroup of finite index in $G$.
\end{enumerate}
\end{cor}

\begin{proof}[Proof of Theorem \ref{central}]
It follows from results of \S \ref{counter} that every very simple
representation enjoys all the properties (i)-(iv). Now suppose
that an absolutely irreducible representation
$$\rho:G \to \Aut_k(V)$$
enjoys the properties (ii)-(iv). It follows from Remark
\ref{proj}(iii) that the $G$-module $V$ does not split.

 Let $R\subset\End_k(V)$ be
a normal subalgebra. Since the $G$-module $V$ is not induced, it
follows from Lemma \ref{induce} and Example \ref{ind} that the
faithful $R$-submodule $V$ is isotypic. This means that
 there exist a simple $R$-module $W$, a positive integer $d$ and an isomorphism
 $$\psi:V \cong W^d$$
of $R$-modules. The following arguments are inspired by another
result of Clifford \cite{Clifford}, ~\cite[Satz 17.5 on p.
567]{Huppert}.

Let us put $$V_1=W, \quad V_2=k^d.$$ The isomorphism $\psi$ gives
rise to the isomorphism of $k$-vector spaces $$V=W^d=W\otimes_{k}
k^d= V_1\otimes_{k} V_2.$$ We have
    $$d \cdot \dim_k(W)=\dim_k(V)$$
Clearly, $\End_R(V)$ is isomorphic to the matrix algebra
$\Mat_d(\End_R(W))$  of size $d$ over $\End_R(W)$.

Let us put
    $$k'=\End_R(W).$$
Since $W$ is simple, $k'$ is a finite-dimensional division algebra
over $k$. Since the Brauer group of $k$ is trivial, $k'$ must be a
field. Clearly, $k'$ is a finite algebraic extension of $k$.

We have

    $$\End_k(V)\supset \End_R(V) = \Mat_d(k')\supset k'.$$
In particular,
$$k \subset k'\subset \End_k(V).$$

Clearly, $\End_R(V) \subset \End_{k}(V)$ is stable under the
adjoint action of $G$. This induces a homomorphism

    $$\alpha: G \to\Aut_{k}(\End_R(V))=\Aut_{k}(\Mat_d(k')),
    \ \alpha(s)(u)=\rho(s)u \rho(s)^{-1} \quad \forall s\in
G, u\in \End_R(V).$$

 Since $k'$ is the center of $\Mat_d(k')$, it is stable under the conjugate action of $G$.
 Thus we get a homomorphism $\chi:G\to\Aut(k'/k)$ such that
$$\chi(s)(a)=\alpha(s)(a)=\rho(s)a \rho(s)^{-1} \quad \forall s\in
G, a\in k'.$$
  I claim that the
 absolute irreducibility of $V$ implies that $k'/k$ is Galois and
 $\chi$ is surjective. Indeed, the inclusion $k'\subset \End_k(V)$
 provides $V$ with a natural structure of $k'$-vector space and it
 is clear that
$$\rho(s)(av)=(\chi(s)(a))v \quad \forall s \in G, a \in k', v
\in V.$$ It follows from Remark \ref{semi} that $k'/k$ is Galois
and $\chi:G \to\Aut(k'/k)=\Gal(k'/k)$ is surjective. Since $V$
does {\sl not} admit a twisted multiplication, $k'=k$.

 This implies
that $\End_R(V) = \Mat_d(k)$ and one may rewrite $\alpha$ as

$$\alpha: G\to \Aut_{k}(\Mat_d(k))=\Aut(\End_{k}(V_2))=
\Aut_{k}(V_2)/k^*=\PGL(V_2).$$

 It follows from the Jacobson density theorem that
 $R=\End_{k}(W)=\End_k(V_1)$.

The adjoint action of $G$ on $R$ gives rise to a homomorphism
$$\beta: G \to \Aut_{k}(\End_{k}(W))=\Aut_{k}(\End_{k}(V_1)) =\Aut_{k}(V_1)/k^*=
\PGL(V_1).$$

  Notice that
$$R=\End_{k}(V_1)=\End_{k}(V_1)\otimes
\I_{V_2}\subset\End_{k}(V_1)\otimes_{k}
\End_{k}(V_2)=\End_{k}(V).$$

Clearly, there exists a central extension $\pi:X
\twoheadrightarrow G$ such that one may lift {\sl projective}
representations $\alpha$ and $\beta$ to linear representations

$$\rho'_2: X \to \Aut_{k}(V_2), \quad \rho_1: X \to
\Aut_{k}(V_1).$$ respectively. For instance, one may take as $X$
the subgroup of $G\times\Aut_{k}(V_1)\times\Aut_{k}(V_2)$ which
consists of all triples $(g, u_1, u_2)$ such that $\alpha(g)$
coincides with the image of $u_2$ in $\Aut_{k}(V_2)/k^*$ and
$\beta(g)$ coincides with the image of $u_1$ in
$\Aut_{k}(V_1)/k^*$. In this case the homomorphisms $\pi,
\rho_1,\rho'_2$ are just the corresponding projection maps
$$(g, u_1, u_2)\mapsto g;\quad (g, u_1, u_2)\mapsto u_1;\quad (g, u_1, u_2)\mapsto
u_2.$$
 Now I am going to check that the tensor product
$\rho_1\otimes\rho'_2$ coincides with the composition
$$\rho\pi:X\twoheadrightarrow G\to\Aut_k(V)$$
up to a twist by a linear character of $X$.

 In order to do that, notice that if $x\in X$ and  $g=\pi(x)\in G$
 then
 the conjugation by $\rho(g)$ in
$\End_{k}(V)=\End_{k}(V_1\otimes_{k}V_2)$ leaves stable
 $R=\End_{k}(V_1)\otimes\I_{V_2}$ and coincides
 on $R$ with the conjugation by $\rho_1(x)\otimes\I_{V_2}$ (by the definition of $\beta$ and $\rho_1$).
 Since the centralizer of $\End_{k}(V_1)\otimes \I_{V_2}$ in
$$\End_{k}(V)=\End_{k}(V_1)\otimes_{k} \End_{k}(V_2)$$
coincides with $\I_{V_1}\otimes \End_{k}(V_2)$, there exists $u\in
\Aut_{k}(V_2)$ such that $$\rho(g)=\rho_1(x)\otimes u.$$ Since the
conjugation by $\rho(g)$ leaves stable the centralizer of $R$,
i.e. $\I_{V_1}\otimes \End_{k}(V_2)$ and coincides on it with the
conjugation by $\I_{V_1}\otimes \rho'_2(x)$ (by the definition of
$\alpha$ and $\rho'_2$), there exists a non-zero constant
$\lambda=\lambda(x) \in k^*$ such that $u=\lambda \rho'_2(x)$.
This implies that for each $x\in X$ there exists a non-zero
constant $\lambda=\lambda(x)$ such that
$$\rho\pi(x)=\rho(g)=\rho_1(x)\otimes u=\lambda\cdot
\rho_1(x)\otimes\rho'_2(x).$$ Since both
$$\rho\pi: X \to \Aut_k(V), \quad \rho_1\otimes\rho'_2: X \to
\Aut_k(V),$$ are group homomorphisms, one may easily check that
the map
$$X \to k^*, \quad x \mapsto \lambda=\lambda(x)$$
is a group homomorphism (linear character). Let us define $\rho_2$
as the {\sl twist}
$$\rho_2: X \to \Aut_k(V),\quad
\rho_2(x)=\lambda(x)\cdot\rho'_2(x)\quad \forall x\in X.$$
Clearly, $\rho_2$ is a linear representation of $X$ and
$$\rho\pi=\rho_1\otimes\rho_2.$$
Since the $G$-module $V$ does not split, either $\dim_k(V_1)=1$ or
$\dim_k(V_2)=1$. If $\dim_k(V_1)=1$ then
$R=\End_k(W)=\End_k(V_1)=k$ consists of scalars. If
$\dim_k(V_2)=1$ then $d=\dim_k(V_2)=1$, i.e., $V=W$ and
$R=\End_k(W)=\End_k(V).$
\end{proof}

\section{Doubly transitive permutation groups}
\label{permute}
 Let $B$ be a finite set consisting of $n \ge 3$
elements. We write $\Perm(B)$ for the group of all permutations of
$B$. A choice of ordering on $B$ gives rise to an isomorphism
$$\Perm(B) \cong \Sn.$$
Let $G$ be a  subgroup of $\Perm(B)$. For each $b \in B$ we write
$G_b$ for the stabilizer of $b$ in $G$; it is a subgroup of $G$.

Let $k$ be a field. We write $k^B$ for the $n$-dimensional
$k$-vector space of maps $h:B \to k$. The space $k^B$ is provided
with a natural action of $\Perm(B)$ defined as follows. Each $s
\in \Perm(B)$ sends a map
 $h:B\to \F$ into  $sh:b \mapsto h(s^{-1}(b))$. The permutation module $k^B$
contains the $\Perm(B)$-stable hyperplane
$$(k^B)^0=
\{h:B\to k\mid\sum_{b\in B}h(b)=0\}$$ and the $\Perm(B)$-invariant
line $k \cdot 1_B$ where $1_B$ is the constant function $1$. The
quotient $k^B/(k^B)^0$ is a trivial $1$-dimensional
$\Perm(B)$-module.

Clearly, $(k^B)^0$ contains $k \cdot 1_B$ if and only if
$\fchar(k)$ divides $n$. If this is {\sl not} the case then there
is a $\Perm(B)$-invariant splitting
$$k^B=(k^B)^0 \oplus k \cdot 1_B.$$
Clearly, $k^B$ and $(\F^B)^0$  carry natural structures of
$G$-modules.

Now, let us consider the case of $k=\F_2$. If  $n$ is even then
let us define the $G$-module
$$Q_B:=(\F_2^B)^0/(\F_2 \cdot 1_B).$$
If $n$ is odd then let us put
$$Q_B:=(\F_2^B)^0.$$

\begin{rem}
\label{d2}
 Clearly, $\dim_{\F_2}(Q_B)=n-1$ if $n$ is odd and
$\dim_{\F_2}(Q_B)=n-2$ if $n$ is even. In both cases
$\dim_{\F_2}(Q_B)\ge 2$. One may easily check that $Q_B$ is a
faithful $G$-module if $n \ne 4$.
\end{rem}

The $G$-module $Q_B$ is called the {\sl heart} over the field
$\F_2$ of the group $G$ acting on the set $B$ \cite{Mortimer}. The
aim of this section is to find out when the $G$-module $Q_B$ is
very simple. It follows from Example \ref{dim22} that if $Q_B$ is
very simple then $\dim_{\F_2}(Q_B)> 2$ and therefore $n \ge 5$.

\begin{rem}
\label{odd}
 Assume that $n$ is odd. Then one may easily check that
$\End_G(Q_B)=\F_2$ if and only if $G$ is doubly transitive. This
implies that if $n$ is odd and the $G$-module $Q_B$ is absolutely
simple then $G$ is doubly transitive. This implies that if $n$ is
odd and the $G$-module $Q_B$ is very simple then $G$ is doubly
transitive.
\end{rem}

\begin{rem}
\label{even}
 Let us assume that $n\ge 5$ is even, the $G$-module $Q_B$ is absolutely simple but $G$ is {\sl not} transitive. Let us present $B$ as
a disjoint union of two non-empty $G$-invariant subsets $B_1$ and
$B_2$. Suppose each $B_i$ contains, at least, $2$ elements.
Without loss of generality we may assume that $\#(B_2)\ge \#(B_1)$
and therefore $\#(B_2)\ge 3$.

There is an embedding of $G$-modules
$$\kappa:\F_2^{B_1} \hookrightarrow (\F_2^B)^0, \quad h \mapsto \kappa(h)$$
defined as follows.
$$\kappa(h)(b_1)=h(b_1) \quad \forall b_1\in B_1;\quad \kappa(h)(b_2)=\sum_{b\in
B_1}h(b)\quad \forall b_2\in B_2.$$

Suppose $1_B \in \kappa(\F_2^{B_1})$. Then  both $B_1$ and $B_2$
consist of odd number of elements and therefore
$$\#(B_1)\ge 3, \quad \#(B_2) \ge 3.$$
Clearly,
$$2\le \#(B_1)-1 =\dim_{\F_2}(\kappa(\F_2^{B_1}))=n-\#(B_2)-1 \le n-4<\dim_{\F_2}(Q_B).$$
Therefore the $G$-module $Q_B$ is {\sl not} simple. Contradiction.

Now suppose $1_B$ does not lie  in $\kappa(\F_2^{B_1})$. Then
$$1\le
\#(B_1)=\dim_{\F_2}(\kappa(\F_2^{B_1}))=\#(B_1)=n-\#(B_2)\le
n-3<\dim_{\F_2}(Q_B).$$ Therefore the $G$-module $Q_B$ is {\sl
not} simple. Contradiction.

This implies that either $B_1$ or $B_2$ is a singleton.

We conclude that if $n$ is even, the $G$-module $Q_B$ is simple
but $G$ is not transitive then $B$ is the disjoint union of two
$G$ orbits of cardinality $n-1$ and $1$ respectively. In other
words, there exists $b\in B$ such that $G=G_b$ and the action of
$G$ on $B\setminus \{b\}$ is transitive. Notice that if we denote
$B\setminus {b}$ by $B'$ then the $G$-modules $Q_B$ and $Q_{B'}$
are isomorphic ~\cite[Remark 2.5 on p. 95]{ZarhinCrelle} (see also
 ~\cite[Hilffsatz 3b]{Klemm}). Applying Remark \ref{odd} to $B'$, we
 conclude that the action of $G$ on $B'$ is doubly transitive.
\end{rem}

\begin{rem}
Suppose $n=2m$ is even, $G$ is transitive but not doubly
transitive. Assume also the $G$-module $Q_B$ is very simple. Then
$n\ge 5$ and  $Q_B$ is absolutely simple. According to ~\cite[Satz
11]{Klemm}, this implies that $m$ is odd and there exists a
subgroup $H\subset G$ of index $2$ such that $B$ can be presented
as a disjoint union of two $H$-invariant subsets $B_1$ and $B_2$
of cardinality $m$.

Since $n \ge 5$, we conclude that $m \ge 3$. There is an embedding
of $H$-modules
$$\kappa:\F_2^{B_1} \hookrightarrow (\F_2^B)^0, \quad h \mapsto \kappa(h)$$
defined as follows.
$$\kappa(h)(b_1)=h(b_1) \quad \forall b_1\in B_1;\quad \kappa(h)(b_2)=\sum_{b\in
B_1}h(b)\quad \forall b_2\in B_2.$$

Clearly, $1_B \in \kappa(\F_2^{B_1})$. We have
$$2\le m-1=\dim_{\F_2}(\kappa(\F_2^{B_1}))=m-1<2m-2=n-2=\dim_{\F_2}(Q_B).$$
Therefore the $H$-module $Q_B$ is {\sl not} simple. Since $H$ is
obviously normal in $G$ and the $G$-module $Q_B$ is very simple,
we conclude that $H$ acts on $Q_B$ via scalars. Since
$\F_2^*=\{1\}$, $H$ acts on $Q_B$ trivially. But this contradicts
to the faithfulness of the $G$-module $Q_B$.

We conclude that if $n\ge 5$ is even, $G$ is transitive and the
$G$-module $Q_B$ is very simple then $G$ must be doubly
transitive.
\end{rem}

To summarize, we arrive to the following conclusion.

\begin{thm}
\label{simple}
 Suppose that $n \ge 3$ is an integer, $B$ is an
$n$-element set, $G \subset \Perm(B)$ is a permutation group.
Suppose that the $G$-module $Q_B$ is very simple. Then $n\ge 5$
and one of the following two conditions holds:
\begin{enumerate}
\item[(i)]
$G$ acts doubly transitively on $B$;
\item[(ii)]
$n$ is even, there exists a $G$-invariant element $b\in B$ and $G$
acts doubly transitively on $B':=B\setminus \{b\}$. In addition,
the $G$-modules $Q_B$ and $Q_{B'}$ are isomorphic.
\end{enumerate}
\end{thm}

\begin{ex}
\label{Lnq} Suppose that there exist a positive integer $m>2$ and
an odd power prime $q$ such that $n=\frac{q^m-1}{q-1}$ and one may
identify $B$ with the $(m-1)$-dimensional projective space
$\P^{m-1}(\F_q)$ over $\F_q$ in such a way that $G$ contains
$\L_m(q)=\PSL_m(\F_q)$. Then the $G$-module $Q_B$ is very simple.

Indeed, in light of Remark \ref{image}(ii), we may assume that
$G=\L_m(q)$; in particular $G$ is a simple non-abelian group
acting doubly transitively on $B=\P^{m-1}(\F_q)$.

Assume that $(m,q) \ne (4,3)$. It follows from a result of
Guralnick \cite{GurTiep} that every nontrivial projective
representation of $G=\L_m(q)$ in characteristic $2$ has dimension
$\ge \dim_{\F_2}(Q_B)$ (see ~\cite[Remark 4.4]{ZarhinMMJ}). It
follows from Example \ref{small} that the $G$-module $Q_B$ is very
simple.

So, we may assume that
 $m=4,q=3$. We have $n=\#(B)=40$ and $\dim_{\F_2}(Q_B)=38$.
It is known \cite{Mortimer} that the $G$-module $Q_B$ is
absolutely simple.

According to the Atlas ~\cite[pp. 68-69]{Atlas}, $G=\L_4(3)$ has
two conjugacy classes of maximal subgroups of index $40$. All
other maximal subgroups have index greater than $40$.
 Therefore all proper subgroups of $G$  have index greater than
 $39>38$ and therefore $Q_B$ is not induced from a representation
 of a proper subgroup.

It follows from the Table on p. 165 of \cite{AtlasB} that all
absolutely irreducible representations
 of $G$ in characteristic $2$ have dimension which is {\sl not}
 a strict divisor of $38$. Applying Corollary \ref{centralF2}, we
 conclude that $Q_B$ is very simple.
\end{ex}

\begin{thm}
\label{double} Suppose that $n\ge 5$ is an integer, $B$ is a set
consisting of $n$ elements. Suppose $G \subset \Perm(B)$ is one of
the known doubly transitive permutation groups (listed in
\cite{Mortimer,DM}). Then the $G$-module $Q_B$ is very simple if
and only if one of the following conditions holds:
\begin{enumerate}
\item[(i)]
$G$ is isomorphic either to the full symmetric group $\SS_n$ or to
the alternating group $\A_n$;
\item[(ii)]
There exist a positive integer $m>2$ and an odd power prime $q$
such that $n=\frac{q^m-1}{q-1}$ and one may identify $B$ with the
$(m-1)$-dimensional projective space $\P^{m-1}(\F_q)$ over $\F_q$
in such a way that $G$ contains $\L_m(q)=\PSL_m(\F_q)$;
\item[(iii)]
$q=n-1$ is a power of $2$ and one may identify $B$ with the
projective line $\P^1(\F_q)$ in such a way that $G$ contains
$\L_2(\F_q)=\PSL_2(\F_q)$;
\item[(iv)]
There exists a positive integer $d\ge 2$ such that
$q:=2^d,n=q^3+1$ and $G$ contains a subgroup isomorphic to the
projective special unitary group $\UU_3(q)=\P\SU(3,\F_{q^2})$;
\item[(v)]
There exists a positive integer $d\ge 2$ such that $q=2^{2d+1},
n=q^2+1$ and $G$ contains a subgroup isomorphic to the Suzuki
group $\Sz(q)$;
\item[(vi)]
$n=11$ and $G$ is isomorphic either to $\L_2(11)$ or to the
Mathieu group $\M_{11}$;
\item[(vii)]
$n=12$ and $G$ is isomorphic either to $\M_{11}$ or to the Mathieu
group $\M_{12}$.
\end{enumerate}
\end{thm}

\begin{proof}
The fact that all the $G$-modules $Q_B$  arised from
\ref{double}(i)-(vii) are very simple was proven in
\cite{ZarhinTexel}(cases (i), (iii), (v), (vi), (vii)),
\cite{ZarhinPAMS}(case (iv)) and in the present paper (Example
\ref{Lnq}: case (ii)). On the other hand,  the paper
\cite{Mortimer} (complemented by \cite{Ivanov}) contains the list
of  doubly transitive $G\subset \Perm(B)$ with absolutely simple
$Q_B$. In addition to the cases \ref{double}(i)-(vii), the
$G$-module $Q_B$ is absolutely simple only if one of the following
conditions holds:
\begin{enumerate}
\item[(a)]
There exists an odd power prime $q$ and a positive integer $d$
such that $n=q^d$ and one may identify $B$ with the affine space
$\F_q^d$ in such a way that $G$ is contained in
$\A\Gamma\L(d,\F_q)$ and contains the group $\F_q^d$ of
translations.  Here $\A\Gamma\L(d,\F_q)$ is the group of
permutations of $\F_q^d$ generated by the group $\AGL(d,\F_q)$ of
affine transformations and the Frobenius automorphism;
\item[(b)]
There exists an odd power prime $q$ such that $n=q+1$ and one may
identify $B$ with the projective line $\P^1(\F_q)$ in such a way
that $G$ becomes a $3$-transitive subgroup of
$\P\Gamma\L(2,\F_q)$. Here $\P\Gamma\L(2,\F_q)$ is the group of
permutations of $\P^1(\F_q)$ generated by $\PGL(2,\F_q)$ and the
Frobenius automorphism.
\end{enumerate}
In the case (a) the group $\F_q^d$ of translations is a proper
normal abelian subgroup of $G$. It follows from Remark
\ref{image}(vi) that $Q_B$ is not very simple.

In the case (b) let us consider the intersection $G'=G\bigcap
\PSL(2,\F_q)$. Clearly, $G'$ is a normal subgroup of $G$. Since
the $\PSL(2,\F_q)$-module $Q_B$ is not absolutely simple
\cite{Mortimer}, the $G'$-module is also not absolutely simple. By
Remark \ref{image}(vi),  $G'$ acts on $Q_B$ by scalars. Since
$\F_2^{*}=\{1\}$ and the $G$-module $Q_B$ is faithful, $G'=\{1\}$.
Since $\PSL(2,\F_q)$ is a subgroup of index $2$ in $\PGL(2,\F_q)$,
the intersection $H:= G\bigcap \PGL(2,\F_q)$ is either a normal
subgroup of order $2$ in $G$ or trivial (one-element subgroup). In
the latter case $G$ is isomorphic to a subgroup of the cyclic
quotient $\P\Gamma\L(2,\F_q)/\PGL(2,\F_q)$ and therefore is
commutative which contradicts the absolute simplicity of the
$G$-module $Q_B$. In the former case, $H$ is an  abelian normal
subgroup of $G$ and it follows from  Remark \ref{image}(vi) that
$Q_B$ is not very simple.
\end{proof}

\section{Applications to hyperelliptic jacobians}
Throughout this section we assume that $K$ is a field of prime
characteristic $p$ different from $2$. We fix its algebraic
closure $K_a$ and write $\Gal(K)$ for the absolute Galois group
$\Aut(K_a/K)$. Let $n \ge 5$ be an integer. Let $f(x)\in K[x]$ be
a polynomial of degree $n$ without multiple roots. We write
$\RR_f$ for the set of roots of $f$. Clearly, $\RR_f$ is a subset
of $K_a$ consisting of $n$ elements. Let $K(\RR_f)\subset K_a$ be
the splitting field of $f$. Clearly, $K(\RR_f)/K$ is a Galois
extension and we write $\Gal(f)$ for its Galois group
$\Gal(K(\RR_f)/K)$. By definition, $\Gal(K(\RR_f)/K)$ permutes
elements of  $\RR_f$; further we identify $\Gal(f)$ with the
corresponding subgroup of $\Perm(\RR_f)$.

\begin{rem}
\label{irr}
 Clearly, $\Gal(f)$ is transitive if and only if the
polynomial $f(x)$ is irreducible. It is also clear that the
following conditions are equivalent:
\begin{enumerate}
\item[(i)]
There exists a root $\alpha \in K$ and an irreducible polynomial
$f_1(x)\in K[x]$ of degree $n-1$ and without multiple roots such
that
$$f(x)=(x-\alpha) f_1(x);$$
\item[(ii)]
There exists a $\Gal(f)$-invariant element $\alpha\in\RR_f$ such
that $\Gal(f)$ acts transitively on $\RR_f\setminus \{\alpha\}$.
\end{enumerate}
\end{rem}

Let $C_f$ be the hyperelliptic curve $y^2=f(x)$. Its genus $g$ is
$\frac{n-1}{2}$ if $n$ is odd and  $\frac{n-2}{2}$ if $n$ is even.
Let $J(C_f)$ be the jacobian of $C_f$; it is a $g$-dimensional
abelian variety defined over $K$. Let $J(C_f)_2$ be the kernel of
multiplication by $2$ in $J(C_f)(K_a)$; it is $2g$-dimensional
$\F_2$-vector space provided with the natural action
$$\Gal(K) \to \Aut_{\F_2}(J(C_f)_2)$$
of $\Gal(K)$. It is well-known (see for instance
\cite{ZarhinTexel}) that the homomorphism $\Gal(K) \to
\Aut_{\F_2}(J(C_f)_2)$ factors through the canonical surjection
$\Gal(K) \twoheadrightarrow \Gal(K(\RR_f)/K)=\Gal(f)$ and the
$\Gal(f)$-modules $J(C)_2$ and $Q_{\RR_f}$ are isomorphic. It
follows easily that the $\Gal(K)$-module $J(C_f)_2$ is very simple
if and only if the $\Gal(f)$-module $Q_{\RR_f}$ is very simple.
Combining Theorem \ref{simple} and Remark \ref{irr}, we obtain the
following statement.

\begin{thm}
\label{ddouble}
 Suppose the $\Gal(K)$-module $J(C_f)_2$ is very
simple. Then one of the following conditions holds:
\begin{enumerate}
\item[(i)]
The polynomial $f(x)\in K[x]$ is irreducible and its Galois group
$\Gal(f)$ acts doubly transitively on $\RR_f$;
\item[(ii)]
$n$ is even, there exists a root $\alpha \in K$  of $f$ and an
irreducible polynomial $f_1(x)\in K[x]$ of degree $n-1$ and
without multiple roots such that
$$f(x)=(x-\alpha) f_1(x).$$
In addition, the Galois group $\Gal(f_1)$ of $f_1$ acts doubly
transitively on $\RR_{f_1}=\RR_{f}\setminus \{\alpha\}$.
%Also, the $\Gal(K)$-modules $J(C_f)_2$ and $J(C_{f_1})_2$ are isomorphic.
%Here $C_{f_1}$ is the hyperelliptic curve $y^2=f_1(x)$.
\end{enumerate}
\end{thm}

\begin{rem}
In the case \ref{ddouble}(ii) the hyperelliptic curves $C_f$ and
$C_{f_1}:y^2=f_1(x)$ are birationally isomorphic over $K$.
\end{rem}

\end{document}